%% file: symprev3.tex
\documentstyle[psfig]{article}

\def\0{{\bf 0}}

\def\bt{\begin{theorem}}
\def\et{\end{theorem}}
\def\bp{\begin{proposition}}
\def\ep{\end{proposition}}
\def\bl{\begin{lemma}}
\def\el{\end{lemma}}
\def\bi{\begin{itemize}}
\def\ei{\end{itemize}}
\def\bd{\begin{description}}
\def\ed{\end{description}}
\def\br{\begin{remark}}
\def\er{\end{remark}}
\def\be{\begin{equation}}
\def\ee{\end{equation}}
\def\bc{\begin{corollary}}
\def\ec{\end{corollary}}
\def\bex{\begin{example}}
\def\eex{\end{example}}

\input mssymb12.tex
\input newmacro.tex

\title{Symplectic Geometry and the Uniqueness of Grauert tubes}

\author{D. Burns\thanks{Research in part supported by NSF grant DMS9408994} \and R. Hind}

\begin{document}

\maketitle

\section{Introduction}

Given a differentiable manifold $M$, it is always possible to
find an almost complex structure $J$ on $TM$ with respect to which the
zero-section is totally-real. In this paper we study a situation in
which there is a
canonical choice for a complex structure with this property, a canonical complexification of $M$.

Suppose that $(M,g)$ is a real-analytic Riemannian manifold of dimension $n$. Identify $M$ with the zero section in $TM$. 
Let
$\rho :TM\to \R$ be the length, with respect to $g$, of tangent
vectors. Then, for $r$ sufficiently small, $T^r M=\{v\in TM | \rho (v) < r\}$
carries a unique complex structure satisfying the following two conditions.

$(i)$ $\rho ^2$ is strictly plurisubharmonic and the corresponding K\"{a}hler
metric restricts to $g$ on $M$.

$(ii)$ $\rho$ is a solution of the homogeneous complex Monge-Amp\`{e}re equation $(dd^c\rho)^n = 0$ on
$T^r M\setminus M$, where $M \subset T^r M$.

This is proven by V. Guillemin and M. Stenzel in \cite{gs} and by L. Lempert and R. Sz\"{o}ke in \cite{ls}. The resulting Stein complex manifolds
will be
called Grauert tubes. Alternatively, one can think of the complex structure
as the unique choice making the leaves of the Riemann foliation on $T^r M \setminus M$ into
holomorphic curves.

It will be important for our purposes that the map
$$\sigma : T^r M \to T^r M, v \mapsto -v$$
is an antiholomorphic involution with respect to these special
complex structures.

Another characterisation of Grauert tubes is given by the following theorem,
taken from \cite{burns}.

\begin{theorem}
Let $X$ be a connected complex manifold and $u$ a smooth, bounded, non-negative,
strictly plurisubharmonic exhaustion function such that
$\sqrt u$ solves
the Monge-Amp\`{e}re equation on $X\backslash \{u=0\}$.

Then $M=\{u=0\}$ is a connected real-analytic submanifold of $X$ and
there exists a biholomorphism $\phi$ from $X$ to a Grauert tube on $T^rM$
with $\rho = \sqrt u\circ \phi$, and $r = \sup_{x \in X} \sqrt{u(x)}$.
\end{theorem}
We call the $r$ in the statement of the theorem the radius of the
Grauert tube.  
The purpose of this paper is to demonstrate that there is a unique way
of associating a Riemannian manifold to a Grauert tube. Specifically, the
remainder of the paper is devoted to proving the following.

\begin{theorem}
Suppose that two Grauert tubes, $X_1$ associated to a Riemannian manifold
$(M_1, g_1)$ and $X_2$ associated to $(M_2, g_2)$, of equal radius
$r$, are biholomorphic
via a map $\phi$. Then $\phi$ maps $M_1 \subset X_1$ to $M_2 \subset X_2$,
restricting to an isometry with respect to $g_1$ and $g_2$.

In particular, the biholomorphisms of a Grauert tube are just the
differentials of the isometries of the underlying Riemannian manifold.
\end{theorem}
The last statement of the theorem is simply the functoriality of the
Grauert tube construction \cite{ls}: every isometry $\phi$ of $M$
extends uniquely by its differential $d\phi_{*}: T^r(M) \rightarrow
T^r(M)$ to a biholomorphic map of the Grauert tube.

A weak form of Theorem 2 was shown in \cite{burns}. Several other
previous partial results are due to S.-J. Kan \cite{kan1, kan2} and
S.-J. Kan and D. Ma \cite{km1, km2}, which employed an interesting variety of
different methods.

The case of Grauert tubes of infinite radius is quite different from
the case treated here. In particular, biholomorphic automorphism
groups for such a Grauert tube can be infinite dimensional, and can
move the submanifold $M$ off itself in the tube. Simple examples of
this phenomenon are shown in \cite{szo2}. Partial uniqueness results
({\em i.e.}, under more restrictive hypotheses) in the  case of
infinite radius appear in \cite{burns} and \cite{szo2}.

Finally, it is a pleasure for us to express our thanks to Yasha
Eliashberg for several very stimulating conversations we have had with
him during the course of this work, from which we have benefited
significantly.

\section{Proof of Theorem 2}

Let $X$ be a Grauert tube of finite radius $r$ associated to a
Riemannian manifold $(M,g)$. Let $Aut_{\Bbb C}(X)$ denote the group of
biholomorphisms of $X$, let $\pm Aut_{\Bbb C} (X)$ be the group of
holomorphic or antiholomorphic automorphisms of $X$,
and $Isom(M) \hookrightarrow Aut_{\Bbb C}(X)$
be the group of isometries of $M$ with respect to the original metric
$g$, acting via differentials on $T^rM$. We note that the index of
$Aut_{\Bbb C}(X)\subset \pm Aut_{\Bbb C}(X)$ is exactly two as $X$
has an antiholomorphic involution.
It follows from Theorem 6.3
of Sz\H{o}ke \cite{szo2} that $Aut_{\Bbb C}(X)$,
and hence $\pm Aut_{\Bbb C}(X)$, is a compact Lie
group.

Now, following \cite{burns}, we let $\tau =\rho ^2$, the length function
squared, and proceed as follows. The function $\tau :X\to [0,r^2)$ is a
proper strictly plurisubharmonic exhaustion function of $X$. Let
$\tilde{\tau}$ be the average of $h^* (\tau)$ over $h\in \pm Aut_{\Bbb C}(X)$
with respect to the Haar measure on $\pm Aut_{\Bbb C}(X)$. Then, as
explained in \cite{burns}, $\tilde{\tau}$ is again a smooth, strictly
plurisubharmonic exhaustion of $X$ which is now invariant under
$\pm Aut_{\Bbb C}(X)$. Therefore choosing $c$ to be slightly less
than $r^2$, the set $\tilde{X} =\{z\in X|\tilde{\tau}\le c\}$ is a
strictly pseudoconvex subset of $X$ which is invariant under
$\pm Aut_{\Bbb C}(X)$. As $\tilde{X}$ lies in a slightly larger
Stein manifold, namely $X$, we can apply a
theorem of S-Y. Cheng and S-T. Yau, see \cite{yau}, to say that $\tilde{X}$ carries a
unique complete K\"{a}hler-Einstein metric of negative Ricci
curvature. 

\vskip 5pt

{\bf Remark} There is a maximal radius $R$, which is not necessarily
$+\infty$, for which a Grauert tube structure can be defined on
$T^R M$. The above discussion is of course only necessary if our $r=R$.
Otherwise a K\"{a}hler-Einstein metric exists on $X$ itself and
in what follows we can just assume $\tilde{X}=X$.

\vskip 5pt

Let $g_{KE}$ be the real (Riemannian) part of the 
K\"{a}hler-Einstein metric
and $\omega_{KE}$ the imaginary part, which is a symplectic form on
$\tilde{X}$. Since the involution $\sigma$ is anti-holomorphic, we have that
$\sigma ^{*} g_{KE} =g_{KE}$ and $\sigma ^* \omega_{KE}
=-\omega_{KE}$.  The zero-section $M\subset X$ is fixed by $\sigma$,
and hence $M$ is a totally geodesic submanifold of $\tilde{X}$ with the metric
$g_{KE}$. Also, $\omega_{KE}$ must be zero when restricted to $M$, in
other words, $M$ is a Lagrangian submanifold of $(\tilde{X},\omega_{KE})$.

These observations lead to the following lemma (see \cite{burns}).

\begin{lemma}
$Aut^0_{\Bbb C}(X) \simeq Isom^0(M)$.
\end{lemma}

Here, $Aut^0_{\Bbb C}(X)$ denotes the identity component of the Lie group
of biholomorphisms of $X$, and $Isom^0(M)$ is the identity component
of the isometries of $M$ with respect to the original metric $g$, acting via differentials on $T^rM$.

{\bf Proof}

It is a result of Chen, Leung and Nagano, see \cite{cln}, that a
closed, maximal-dimensional, totally-real, minimal submanifold $M$ of
a K\"{a}hler manifold $\tilde{X}$ of negative Ricci curvature is strictly stable.
In our situation this implies in particular that $Aut^0_{\Bbb C}(X)$ must
preserve $M$ as all biholomorphisms act as isometries of $g_{KE}$
when restricted to $\tilde{X}$.
Such biholomorphisms preserving $M$ must act as isometries of $g$
when restricted to $M$, see \cite{szo2}, and it has already been
remarked that the differentials of isometries of $g$ are biholomorphisms
of $X$.

\vskip 5pt

A few comments will be needed on the symplectic manifolds
$(\tilde{X},\omega_{KE})$.

\vskip 5pt

{\bf Remark} We observe that for the remainder of this paper the symplectic
form $\omega_{KE}$ on $\tilde{X}$ or $X$ could in fact be replaced by any
other symplectic form which is invariant under $\pm Aut_{\Bbb C}(X)$ and,
like $\omega_{KE}$, satisfies the conclusions of the forthcoming Lemma $5$.
Such symplectic forms can be constructed directly from the canonical form
on the cotangent bundle of $M$ without appealing to Cheng and Yau's
analysis. We thank the referee for this remark.

\vskip 5pt

Let $h:[0,r) \to [0,\infty)$ be a smooth, convex and strictly
increasing proper map satisfying $h(x)=x$ for $x$ near $0$.

\begin{lemma}
$(\tilde{X},\omega_{KE})$ is symplectomorphic via a symplectomorphism
fixing $M$ to $(T^r M, -dd^cf)$, where $f=h \circ \rho^2$.
\end{lemma}

{\bf Proof}

The K\"{a}hler-Einstein metric can also be written in the form
$\omega_{KE}=-dd^cf'$ for some strictly plurisubharmonic exhaustion function $f'$ on
$\tilde{X}$ (\cite{yau}). Now the existence of a symplectomorphism follows from
Theorem 1.4.A in \cite{elg}.
Since $\sigma$ is anti-holomorphic, we have $\sigma^{*} \circ dd^c = - dd^c
\circ \sigma^{*}$, and thus replacing $f'$, if necessary, by $F =
\frac{1}{2}(f' + \sigma^{*}f')$ we can assume without loss of
generality that the exhaustion function
$f'$ is $\sigma$-invariant. Hence the whole
construction of \cite{elg} can be carried out in a $\sigma$-invariant
fashion. In particular the resulting symplectomorphism must commute
with $\sigma$ and hence fix $M$, the fixed point set of $\sigma$.

\vskip 5pt
 
The symplectic manifold $(T^r M, \omega =-dd^cf)$ has, in
the terminology of \cite{elg}, a contracting vector field $v$
given by the negative gradient of $f$ with respect to the metric
$G(x,y)=\omega(x,Jy)$, where $J$ is our complex structure.
Equivalently, $v$ satisfies
$v \rfloor \omega =d^cf$.
We note that $v$ vanishes only along $M$.

Furthermore, 
$${\cal L}_v \omega =d(v \rfloor \omega) =-\omega$$
and
$${\cal L}_v d^cf=-v \rfloor \omega + d(v \rfloor d^cf)=-d^cf.$$

In particular, if $(v_t)$ denotes the
$1$-parameter group of diffeomorphisms generated by $v$ we have
\begin{equation}
v_t^* d^cf = e^{-t} d^cf. \label{exact}
\end{equation}

We can now in fact observe the following.

\begin{lemma}
$(\tilde{X},\omega_{KE})$ is symplectomorphic via a symplectomorphism $\psi$
fixing $M$ to $(T^* M,d(pdq))$ where $pdq$ denotes the canonical
Liouville $1$-form on the cotangent bundle.
\end{lemma}

{\bf Proof}

It suffices to show that $(T^r M, \omega =-dd^cf)$ is symplectomorphic
to $(TM, d\lambda)$ where $\mu$ is the pull-back of $pdq$ to $TM$
using the isomorphism given by the metric $g$. The conditions on a
Grauert tube imply that $\mu =-d^cf$ in the neighbourhood
of $M$ where $h(x)=x$.
Recall that $f=h\circ \rho^2$ where  $h:[0,r) \to [0,\infty)$ is any smooth, convex and strictly
increasing proper map satisfying $h(x)=x$ for $x$ near $0$.

Now,
$$d^cf=h'(\rho^2)d^c\rho^2$$
and
$$\omega=-dd^cf=-h''(\rho^2)d(\rho^2)\wedge d^c\rho^2-h'(\rho^2)dd^c\rho^2.$$

We choose $h$ to grow sufficiently fast that $h''\gg h'$ for $x$ close
to $r$.
Then if $v$ is the corresponding contracting vector field we see that
$d(\rho^2)(v)$ must be of the order of $\frac{h'}{h''}$ near the
boundary of $T^r M$. Hence by choosing a suitable $h$ we are able to
assume that the vector field $-v$ is complete.

We can define an expanding vector field $w$ on $TM$ with respect to
$d \mu$ by $w \rfloor d\mu = \mu$. Let $(w_t)$ denote the
corresponding $1$-parameter group.

In a neighbourhood $U$ of $M$ we have $w=-v$.

We now define the map $\psi$ of $T^r M$ to $TM$ as follows.
Given $x \in TM$, choose $t$ sufficiently large that
$v_t (x)\in U$. Then set $\psi (x)=w_t \circ v_t (x)$.
This is clearly well-defined (that is, independent of $t$),
is a symplectomorphism (which is surjective since $-v$ is complete),
and is the identity near $M$.

\vskip 5pt

We can now prove Theorem $2$.
Suppose that we have a biholomorphism $\phi$ between two Grauert
tubes $X_1 =T^{r_1}M_1$ and $X_2 =T^{r_2}M_2$.
We will denote $\phi(M_1)$ also by $M_1$, and construct $\tilde{X_2}$
as above. There is an antiholomorphic involution $\sigma_1$ of $X_2$
which preserves $M_1$, namely the push-forward of the involution on $X_1$.
As $\sigma_1$ also preserves $\tilde{X_2}$, we see that $M_1$,
alongside $M_2$, must
be a Lagrangian submanifold of $(\tilde{X_2},\omega_{KE})$.
Combining with the results of \cite{szo1} and
\cite{burns},
to prove Theorem $2$ it will suffice for us to show that as submanifolds
of $X_2$ we have
$M_1 = M_2$.

Let $\lambda$ and $v$ denote the pull-backs via $\psi$ of the primitive $1$-form
$pdq$ and the corresponding contracting vector field respectively
to $\tilde{X_2}$. Also, let $\omega = d\lambda = \omega_{KE}$. 

There are two cases to consider.

\vskip 5pt

{\bf Case $1$} $M_1 \subset \tilde{X_2}$ is not an exact Lagrangian submanifold.
That is, the form $\lambda |_{M_1}$, which is closed because $M_1$ is
Lagrangian, is not exact.

In this case, there is a closed loop $\gamma_1 \subset M_1$ along which
$\lambda$ has a non-zero integral. Project $\gamma_1$ onto $M_2$ by the
bundle projection of $X_2=T^{r_2} M_2$ and connect each point of $\gamma_1$
with its projection by a segment in the corresponding fiber to construct (the
image of) a cylinder $C$ in $\tilde{X_2}$ whose boundary components are
$\gamma_1$ and the projection, say $\gamma_2$, of $\gamma_1$ in $M_2$. As $\lambda |_{M_2} =0$,
we have
$$\int_{C} \omega = \int_{\gamma_1} \lambda \ne 0$$
by Stokes' Theorem.

We recall now the antiholomorphic involution $\sigma$
of a Grauert tube about the underlying Riemannian manifold.
In our case we have antiholomorphic involutions
$\sigma_1$ about $M_1$ and $\sigma_2$ about $M_2$.

In fact there is a sequence of (Lagrangian) submanifolds
of $\tilde{X_2}$ about which there are antiholomorphic involutions.
We define (unfortunately) $N_1=M_2$, $N_2=M_1$, $N_3$ to be
the reflection of $N_1$ in $N_2$, and in general
$N_k$ to be the reflection of $N_{k-2}$ in $N_{k-1}$.
Each $N_k$ can be written in the form
$\sigma_1^{\epsilon} (\sigma_2 \sigma_1)^{n(k)}M_2$ for $k$ odd
and $\sigma_2^{\epsilon} (\sigma_1 \sigma_2)^{n(k)}M_1$ for $k$ even.
In these formulas, $\epsilon$ is either $0$ or $1$ depending on $k(mod 4)$
and similarly $n(k)$ is some exponent depending on $k$.

As already noted above, it follows from Theorem 6.3 of Sz\H{o}ke
\cite{szo2} that $Aut_{\Bbb C}(X)$ is a compact Lie group for any
Grauert tube $X$ of finite radius. Combining this with the fact that
the identity component of $Aut_{\Bbb C}(X_2)$ must preserve $M_2$ we
deduce that for some large (odd) value of $k$, say $l$, $N_l=M_2$.

Let $C_1=C$, $C_2$ be the reflection of $C_1$ in $N_2$, and in
general $C_k$ be the reflection of $C_{k-1}$ in $N_k$.
These cylinders can be joined together end-to-end for $1\le k \le l-1$
to obtain one long cylinder $\tilde{C}$ with both boundaries
in $M_2$. Now, each of the reflections $\sigma$ is an orientation
reversing map of $\tilde{C}$  wherever it is defined, but also
satisfies $\sigma^* \omega = -\omega$. Therefore, we have that  
$$\int_{\tilde{C}} \omega = (l-1)\int_{C} \omega \ne 0.$$

However, this is a contradiction to Stokes' Theorem.

\vskip 5pt

{\bf Case $2$} $M_1 \subset \tilde{X_2}$ is an exact Lagrangian submanifold,
that is, $\lambda |_{M_1}$ is exact.

First we set-up the Lagrangian submanifolds $N_k$ exactly as in Case $1$.
Again assume that $N_l=M_2$.

We now study the isotopy $L_t$ for $t\ge 0$ of $M_1$ given by
$L_t =v_t(M_1)$. Recall that $(v_t)$ is the $1$-parameter group
generated by the contracting vector field $v$.
By equation \ref{exact}, this is an exact Hamiltonian isotopy,
see for instance \cite{chap}.

Suppose that $M_1$ does not map into $M_2$. Then we can find a
point $p\in M_1 \backslash M_2$.

Let $Z= \R +i(0,1) \subset \C$. It is proven by H. Hofer in \cite{hof}
that for any fixed $t$ there exists a holomorphic map $u:Z \to X_2$ with
a continuous extension to the boundary such that
$u(i)=p$, $u(\R)\subset L_t$, $u(\R +i)\subset M_1$ and
$\int_{Z} {u^* \omega} < \infty$.

As $\tilde{X_2}$ is exhausted by pseudoconvex domains, all holomorphic strips
with boundary on $M_1 \cup L_t$, for any $t$, must lie in a
fixed compact region. This observation is sufficient to be able to
apply Hofer's results which are stated for Lagrangian submanifolds
of a compact symplectic manifold. Also we will let $m$ be an
upper bound for the norm of $\lambda$ with respect to the K\"{a}hler
metric on this compact set, which, without loss of generality, includes
all of the $N_k$ for $1\le k \le l$.

In the situation when $M_1$ and $L_t$ do not intersect transversally,
the behaviour of such maps $u(x+iy)$ as $x\to \pm \infty$ may be
difficult to describe precisely. However the following is true
(see \cite{hof} again).

For any $\delta >0$ there exist numbers $R_1, R_2 \in \R$ with
$-R_1 > \delta^{-1}$, $R_2 > \delta^{-1}$ and
$$|\int_{Z} {u^* \omega} - \int_{[R_1,R_2]+i[0,1]} {u^* \omega} | <\delta.$$
Furthermore, for suitable $R_1$, $R_2$ and $j=1,2$, the length of $t\mapsto u(R_j +it)$ is less
than $\delta$. To see this, note that
$$\int_{Z} {(|u_x|^2 +|u_y|^2)dsdt} = 2\int_{Z} {u^* \omega} < \infty.$$

Now,
\begin{eqnarray}
|\int_{Z} {u^* \omega}| & < & \delta + |\int_{[R_1,R_2]+i[0,1]} {u^* \omega}| \nonumber \\
                        & < & (1+2m)\delta + |\int_{[R_1,R_2]} {u^* \lambda} | + |\int_{[R_1,R_2]+i} {u^* \lambda} | \nonumber \\
                        & = & (1+2m)\delta + e^{-t}|\int_{[R_1,R_2]} {u^* (v_t^{-1})^* \lambda} | + |\int_{[R_1,R_2]+i} {u^* \lambda} |. \nonumber
\end{eqnarray}

But the integral of $\lambda$ along any path in $M_1$ is uniformly bounded
(as $\lambda$ is exact it depends only on the endpoints).
Hence, $|\int_{Z} {u^* \omega} |<C$, where $C$ is a uniform constant
independent of $t$.

Let now $u_t$ be the holomorphic map corresponding to the Lagrangian $L_t$,
and for any holomorphic map $f$ from an open subset $E\subset \C$ into
$X_2$ write
$$\p f=\{q\in X_2 | q=\lim_{n\to \infty}f(z_n),z_n \in\E, z_n \to \p E \cup \{\infty\}\}.$$

\begin{lemma}
There exists a sequence $t_j \to \infty$
such that the maps $u_{t_j}$ converge uniformly on compact sets to
a non-constant holomorphic map $u_{\infty}:Z\to X_2$ with $\p u_{\infty} \subset M_1 \cup M_2$.
\end{lemma}

{\bf Proof}

For each $t$, the holomorphic map $u_t$ can be Schwarz reflected in
$M_1$ using $\sigma_1$. After a reparameterization of the resulting
maps such that they are now defined on the open disk $\Delta$ we have
a sequence of holomorphic maps $\tilde{u}_t :\Delta\to X_2$ satisfying
$|\int_{\Delta} {\tilde{u}_t^* \omega} |<2C$, $\tilde{u}_t(0)=p$ and
$\p \tilde{u}_t \subset L_t \cup \sigma_1(L_t)$. Since the
$\tilde{u}_t$ all have their image in a bounded open set in a Stein
manifold, Montel's theorem implies that after taking a subsequence
$t_j$ the $\tilde{u}_{t_j}$ converge uniformly on compact sets to
another holomorphic map $\tilde{u}_{\infty}:\Delta \to X_2$ with
$\tilde{u}_{\infty}(0)=p$. As the $\tilde{u}_t$ all have area
uniformly bounded by $2C$, a result of F. Labourie, see \cite{lab},
implies that $\p \tilde{u}_{\infty}$ is contained in the Hausdorff
limit of the sets $\p \tilde{u}_{t_j}$. This Hausdorff limit is
contained in $M_2 \cup \sigma_1 (M_2)$ which is disjoint from
$p$. Hence the map $\tilde{u}_{\infty}$ is non-constant. This is
enough to establish the lemma, letting $u_{\infty}$ be a suitable
reparameterization of $\tilde{u}_{\infty}$.

\vskip 5pt 

Let now $K= \int_{Z} {u_{\infty}^* \omega}$.
As $u_{\infty}$ is holomorphic and non-constant, $K>0$.

Now choose $\delta$ to be much less than $K$ and find $R_1$ and $R_2$ as
before. As the length of $u_{\infty}(R_j +i(0,1))$ is less than $\delta$, for
$j=1,2$ we can continuously extend $u_{\infty}$ to $R_j +i[0,1]$.

We claim that the map $u_{\infty}$ also extends continuously to $(R_1,
R_2)$, and maps $(R_1, R_2)$ into $M_2 = N_1$, and therefore has a
holomorphic reflection by means of $\sigma_2$ across $(R_1, R_2)$. We
first prepare the target manifold $X_2$.

$X_2$ is a Stein manifold, and so may be embedded properly in ${\Bbb
C}^N$, for some $N$. Let $f_1,...,f_N$ be the component functions of
this embedding, and set $f^{*}_j(z) = \overline{f_j(\sigma_2(z))}, j =
1,...,N$. Consider the embedding $F: X_2 \rightarrow {\Bbb
C}^{2N}$ with components $F_{2j-1} = \frac{1}{2}(f_j + f^{*}_j),
F_{2j} = \frac{1}{2i}(f_j - f^{*}_j), j=1,...,N.$ Note that this
embedding has the property that $F(\sigma_2(z)) = \overline{F(z)}$,
for all $z \in X_2$, and so $F(M_2) = F(X_2) \cap {\Bbb R}^{2N}$. To
prove the map $u_{\infty}: (R_1, R_2) + i(0,1) \rightarrow X_2$
extends by reflection, it suffices to show each of $F_k \circ
u_{\infty}$ extends by reflection to $(R_1, R_2) + i(-1,1)$, and to do
this it suffices, by classical Schwarz reflection, to show that $v_k =
\Im(F_k \circ u_{\infty})$ is continuous up to the boundary arc $(R_1,
R_2)$ and equals $0$ there.  

Now each of the harmonic functions $v_{t_j,k} = \Im(F_k \circ
u_{t_j})$ is continuous up to the boundary of $[R_1, R_2] + i[0,1]$,
and hence can be written as a Poisson integral of its boundary values
there. Note that the functions $v_{t_j,k}$ are uniformly bounded on
$[R_1, R_2] + i[0,1]$. Furthermore, their boundary values converge to
the function $w_{\infty,k}$ which is identically $0$ along $[R_1,
R_2]$ (since $\max_{L_{t_j}} \mid \Im(F_k) \mid $ converges to $0$ as
$t_j \rightarrow \infty$) and which equals $\Im(F_k \circ u_{\infty})$
on the rest of the boundary, by the interior convergence of the
sequence $u_{t_j}$. Hence, by the bounded convergence theorem, we have
that $v_k$ on $(R_1, R_2) + i(0,1)$ is the Poisson integral of
$w_{\infty, k}$. By standard properties of Poisson integrals, $v_k$ is
continuous up to the boundary along the arc $(R_1, R_2)$ in $[R_1,
R_2] + i[0,1]$, completing the proof of the claim.

We now can apply the reflections successively across the $N_k$ exactly
as in Case $1$ to the mapping $u_{\infty}$. Again as in Case $1$ we
stop reflecting when some $N_k=M_2$, so the resulting mapping, which
we still call $u_{\infty}$ is defined and holomorphic on some
rectangle ${\cal R} = [R_1, R_2] + i[0,k]$, and $u_{\infty}([R_1,
R_2]) \subset M_2$ and $u_{\infty}([R_1, R_2] + i \, k) \subset M_2 =
N_1=N_l$.

But $$(l-1)K = \int_{\cal R} u^{*}_{\infty} \omega = \int_{\p {\cal R}}
u^{*}_{\infty} \lambda \le 2m\delta (l-1)$$ giving a contradiction as
required for $\delta$ sufficiently small.

\end{document}

%% file: newmacro.tex
\newcommand{\R}{{\Bbb R}}
\newcommand{\C}{{\Bbb C}}

\newcommand{\E}{{\rm E}}

\def\NABLA#1{{\mathop{\nabla\kern-.5ex\lower1ex\hbox{$#1$}}}}
\def\Nabla#1{\nabla\kern-.5ex{}_#1}

\newcommand{\p}{{\partial}}

\newtheorem{theorem}{Theorem} 
\newtheorem{corollary}[theorem]{Corollary}
\newtheorem{lemma}[theorem]{Lemma}
\newtheorem{proposition}[theorem]{Proposition}

\newtheorem{remark}[theorem]{Remark}
\newtheorem{example}[theorem]{Example}